\documentclass[36pt,reqno]{amsart}

\newcommand\PP{\mathbb P} 
\newcommand{\T}{s_C}  
\newcommand{\W}{\alpha}
 
\renewcommand\O{\mathcal O}

\makeatletter \@addtoreset{equation}{section} \makeatother

\newtheorem{thm}[equation]{Theorem} \newtheorem{lem}[equation]{Lemma}
\newtheorem{cor}[equation]{Corollary}
\newtheorem*{cor*}{Corollary}

\theoremstyle{definition} \newtheorem{defn}[equation]{Definition}
\newtheorem{rmk}[equation]{Remark}

\title[Unstable products of smooth curves]{Unstable products of smooth
  curves} 

\author{J. Ross} 

\begin{document}
\bibliographystyle{plain}

\begin{abstract}
  We give examples of smooth manifolds with negative first Chern class
  which are slope unstable with respect to certain polarisations, and
  so have K\"ahler classes that do not admit any constant scalar
  curvature K\"ahler metrics.  These manifolds also have unstable
  Hilbert and Chow points. We compare this to the work of
  Song-Weinkove on the J-flow.
\end{abstract}

\maketitle

\section{Introduction}

A central problem in K\"ahler geometry is finding necessary and
sufficient conditions for a K\"ahler class on a complex manifold $X$
to admit a constant scalar curvature K\"ahler (cscK) metric.  It is
well known that when the first Chern class of $X$ is negative, there
exists a K\"ahler-Einstein (and hence cscK) metric in $-c_1(X)$
\cite{aub76:equations_monge_ampere,yau78:ricci_kaehler_monge_ampere}.
Thus by a deformation argument due to LeBrun-Simanca
\cite{lebrun_simanca94:_extrem_kaehl} there exists a cscK metric in
each class in an open neighborhood around $-c_1(X)$.

It has been an element of folklore that if $c_1(X)$ is negative then
every K\"ahler class on $X$ admits a cscK metric.  We prove in this paper that
this is not the case by showing that this fails on the product of
certain non-generic smooth curves.  

The cscK problem is related to the stability of $X$.  A conjecture of Yau
states that a rational K\"ahler class $\Omega$ should admit a cscK
metric if and only if the pair $(X,\Omega)$ is ``stable'' in the sense
of geometric invariant theory \cite{yau93:open_problems_in_geometry}.  The
precise definition of stability, called K-stability, was introduced by
Tian
\cite{tian94:k_energy_on_hypersurfaces_and_stability,tian97:kaehler_einstein_metrics_positive}
and expanded by Donaldson \cite
{don02:scalar_curvature_and_stability_of_toric_varieties}.  One
direction of the conjecture has essentially be proved: the
existence of a cscK metric implies stability (see
\eqref{thm:donaldsoncscK}).

To show that certain K\"ahler classes do not admit cscK metrics we
shall use an obstruction for K-stability (and hence for cscK metrics)
of Thomas and the author called slope stability
\cite{rossthomas:_obstructionto_cscK,rossthomas:_hilbert_mumford}.
The main results of this paper are the following: \vspace{.2cm}

{\bf Theorem \ref{cor:main}. }{\it For $g\ge 5$ there exist
  smooth curves $C$ of genus $g$ such that $X=C\times C$ is not slope
  semistable with respect to certain polarisations.  Thus there are
  K\"ahler classes on $X$ that do not admit K\"ahler metrics of
  constant scalar curvature.} \vspace{.2cm}

Moreover this gives the first example of a manifold with negative first
Chern class whose Hilbert and Chow points are unstable in the sense of
geometric invariant theory.  \vspace{.2cm}

{\bf Theorem \ref{thm:unstablehilbertchow}. }{\it For $g\ge 5$ there
  exist smooth curves $C$ of genus $g$ such that $C\times C$ is not
  asymptotically Hilbert semistable (resp.\ not asymptotically Chow
  semistable) with respect to certain polarisations.}  \vspace{.2cm}


More specifically these results holds if $C$ admits a simple branched
cover to $\mathbb P^1$ of degree $d$ with $2\le d-1<\sqrt g$\,
(\ref{thm:maintheorem}).  Such curves are non-generic, since a
general curve of genus $g$ admits a branched cover to $\mathbb P^1$ of
degree
$$d_0=\Big[\frac{g+1}{2}\Big]+1$$ and none of degree $d<d_0$.
(\cite{grif78:principles_algebraic_geometry} p.\ 261).

Given a simple branched cover $\pi\colon C\to \mathbb P^1$ of degree
$d$, consider the fibre product $C\times_{\pi} C\subset X$.  This
contains the diagonal $\Delta$ and we let $Z$ be the residual divisor
(i.e.\ $Z+\Delta=C\times_{\pi} C$). We will show that if $2\le
d-1<\sqrt g$ then $Z$ has negative self-intersection, and the
``slope'' of $Z$ as defined in Section \ref{sec:slope} is less than
that of $X$, which proves instability.

As suggested by Chen, such an example leads one to try and better
understand the obstruction to finding cscK metrics in these classes.
In particular one technique studied by Chen
\cite{chen00:_mabuchilowerbound}, Donaldson
\cite{donaldson99:_momentmapsanddiffeomorphisms}, Weinkove
\cite{weinkove:_jflowmabuchi} and Song-Weinkove
\cite{song_winkove:_jflowmabuchi} is that of the J-flow.  We comment
in Section \ref{sec:mabuchi} as to what some of these results of
\cite{song_winkove:_jflowmabuchi} say when applied to products of
curves.  \vspace{.1cm}

{\bf Acknowledgments} I would like to thank Richard Thomas for
comments on the first draft of this paper, as well as continued
encouragement.  I also thank Xiuxiong Chen, Ian Morrison, Sean Paul,
Michael Thaddeus and Ben Weinkove for useful conversations.
Lazarsfeld's book \cite{lazarsfeld04:_positI} has been extremely
helpful.  \vspace{.1cm}

{\bf Notation } By a polarisation on a variety $X$ we mean a choice of
ample divisor which we usually denote by $L$ and will write $c_1(L)$ for
the first Chern class of the associated line bundle. A $\mathbb
Q$-divisor is a formal sum of divisors with rational coefficients, and
an ample $\mathbb Q$-divisor is one which can be written as a sum of
ample divisors (again with rational coefficients).  The space of
$\mathbb Q$-divisors modulo numerical equivalence is denoted $N^1(X)_\mathbb
Q$.  Abusing notation we will often not distinguish between a $\mathbb
Q$-divisor and its class in $N^1(X)_\mathbb Q$.

\section{Slope stability for  Varieties}\label{sec:slope}

The original link between K-stability and cscK metrics is due to Tian
\cite{tian94:k_energy_on_hypersurfaces_and_stability,
  tian97:kaehler_einstein_metrics_positive}.  For slope stability we
require the definition of K-stability used by Donaldson
\cite{don02:scalar_curvature_and_stability_of_toric_varieties} (the
relation between the two definitions can be found in
\cite{paultian04:_analysisofgeometricstability}).

\begin{thm}\label{thm:donaldsoncscK}
  Fix a polarised manifold $(X,L)$.  If there exists a constant scalar
  curvature K\"ahler metric in $c_1(L)$ then $(X,L)$ is K-semistable.
\end{thm}
\begin{proof}
  This is proved by in \cite{don:lower_calabi}.  Alternatively one can
  use the existence of a cscK metric to show that the Mabuchi
  functional is bounded from below
  \cite{donaldson05:_scalar_ii,chen00:_mabuchilowerbound} which in
  turn implies K-semistability
  \cite{paultian04:_analysisofgeometricstability}.
\end{proof}

The notion of slope stability for polarised varieties was introduced
in \cite{rossthomas:_obstructionto_cscK, rossthomas:_hilbert_mumford}
as a necessary condition for K-stability.  The general idea is that a
non-generic subscheme of a polarised variety $(X,L)$ with certain
numerical properties will have a slope that is too small, and this
forces $(X,L)$ to be unstable.

Let $Z$ be a subscheme of a polarised variety $(X,L)$ and let
$\pi\colon \widehat X\to X$ be the blowup of $X$ along $Z$ with
exceptional divisor $E$.  For sufficiently small positive $c$ the divisor
$\pi^* L-cE$ is ample; so we can define the Seshadri constant of $Z$ as
$$\epsilon(Z,L) = \text{sup}\{ c\in \mathbb Q\colon \pi^* L-cE \text{ is ample}\}.$$  

Write the Hilbert polynomial of $L$ as $\chi(kL) = a_0k^n +
a_1k^{n-1} + \cdots$ where $n=\dim X$.  For fixed $x\in\mathbb Q$ we
define $a_i(x)$ by
$$\chi(k(\pi^*L-xE)) = a_0(x)k^n + a_1(x) k^{n-1} + \cdots \qquad \text{for all } kx\in\mathbb N.$$
As $\chi(k\pi^*L-rE)$ is a polynomial in two variables of total
degree at most $n$, we have that $a_i(x)$ is a polynomial and so
extends to all real $x$.  We let $\tilde{a}_i(x)= a_i -a_i(x)$ and for
$0<c\le \epsilon(Z,L)$ the slope of $X$ and $Z$ is defined to be
(\cite{rossthomas:_obstructionto_cscK} (3.1,3.14)),
\begin{eqnarray*}
  \mu(X,L) &=& \frac{a_1}{a_0},\\
  \mu_c(\O_Z,L) &=& \frac{\int_0^c \tilde{a}_1(x) + \frac{\tilde{a}'_0(x)}{2}dx }{\int_0^c \tilde{a}_0(x) dx},
\end{eqnarray*}
which are both finite (ibid. (4.21)).  

\begin{defn} 
  We say that $(X,L)$ is {\bf slope semistable with respect to a
    subscheme} $Z$ if
    $$\mu(X,L)\le \mu_c(\O_Z,L) \quad \text{ for all } 0<c\le
  \epsilon(Z,L).$$
  We say $(X,L)$ is {\bf slope semistable} if it is
  slope semistable with respect to all subschemes.
\end{defn}

Since the property of being slope semistable is invariant under
replacing $L$ by some power (ibid. (3.10)), we extend the notion of
slope semistability to ample $\mathbb Q$-divisors.

The definition of the slopes are made so that if $(X,L)$ is not slope
semistable with respect to $Z$ then the degeneration to the normal
cone of $Z$ prevents $(X,L)$ being K-semistable:

\begin{thm}\label{thm:slopeimpliesKstability}
  If $(X,L)$ is not slope semistable then it is not K-semistable.
\end{thm}
\begin{proof}
  See \cite{rossthomas:_hilbert_mumford} (4.18);  when
  $X$ and $Z$ are smooth this is proved in
  \cite{rossthomas:_obstructionto_cscK} (4.2).
\end{proof}

\subsection*{Slope stability for smooth surfaces}

In this paper we will only consider the case that $X$ is a smooth
surface and $Z$ is a curve.  Then the blowup of $X$ along $Z$ is just
$X$ itself, so
$$\epsilon(Z,L) = \sup\{ c\in \mathbb Q \colon L-cZ \text{ is ample} \}.$$

Letting $K$ be the canonical divisor of $X$, a simple application of
the Riemann-Roch theorem to calculate $a_0(x)$ and $a_1(x)$ yields
(\cite{rossthomas:_obstructionto_cscK} (5.4)),
\begin{eqnarray}\label{eq:slopeonsurfaces}
  \mu(X,L) &=& -\frac{K.L}{L^2},\\
  \mu_c(\O_Z,L) &=& \frac{3(2L.Z-c(K.Z+Z^2))}{2c(3L.Z-cZ^2)}. \nonumber
\end{eqnarray} 
Notice that in this case $\epsilon(Z,L)$ as well as the slopes
$\mu_c(\O_Z,L)$ and $\mu(X,L)$ depend only on the class of $L$ and $Z$
modulo numerical equivalence.  We extend the equations
\eqref{eq:slopeonsurfaces} to any class $L$ in $N^1(X)_\mathbb Q$
which is not necessarily ample (in which case they may no longer be
finite).

\subsection*{Hilbert and Chow stability}

Slope stability also gives an obstruction to the classical notions of
stability for projective varieties.  For $r\gg 0$ consider the
embedding of $X$ via the linear series $|rL|$ into $\PP^{N(r)}$.  Up
to change of change of coordinates, this determines a point
$\text{Hilb}(X,L^r)$ in the Hilbert scheme of $\PP^{N(r)}$ (resp.\ a
point $\text{Chow}(X,L^r)$ in the Chow variety).  The action of the
automorphism group of $\PP^{N(r)}$ induces a linearised action on the
Hilbert scheme (resp.\ Chow variety), and one can apply the notions of
geometric invariant theory to these spaces (see
\cite{mum77:stability_projective_varieties}).
\begin{defn}
  We say that $(X,L)$ is asymptotically Hilbert (resp.\ Chow)
  semi\-stable if for $r\gg 0$ the point $\text{Hilb}(X,L^r)$ (resp.
  $\text{Chow}(X,L^r)$) is semistable in the sense of geometric
  invariant theory.
\end{defn}

\begin{thm}\label{thm:slopeimplieshilbertchowstability}
  If $(X,L)$ is not slope semistable then it is neither asymptotically
  Hilbert nor asymptotically Chow semistable.
\end{thm}
\begin{proof}
  The follows from \eqref{thm:slopeimpliesKstability} as
  asymptotic Hilbert (resp.\ Chow) semistability implies
  K-semistability (\cite{rossthomas:_hilbert_mumford} (4.18)).
\end{proof}

\section{Unstable products of Curves}

We start with some standard material on the ample cone of products of
curves, all of which can be found in \cite{lazarsfeld04:_positI}.  Fix
a smooth curve $C$ of genus $g\ge 2$ and let $X=C\times C$.  If
$\pi_i$ is the projection onto the $i$-th factor, and $p$ is a fixed
point in $C$ then the class $f_i$ of the fibre $\pi_i^{-1}(p)$ in
$N^1(X)_\mathbb Q$ is independent of $p$.  The class of the canonical
divisor of $X$ is $K=(2g-2)(f_1+f_2)$ which is ample, so $X$ has
negative first Chern class.  Letting $\delta$ be the class of the
diagonal we see that $f_i^2=0$, $f_1.f_2=1$, $f_i.\delta=1$ and
$\delta^2=2-2g$.

Let $f=f_1+f_2$ and for convenience make the change of variables
$\delta'=\delta-f$.  Then we have the following intersection numbers on
$X$:
$$ f^2 =2, \quad \delta'.f=0,\,\,\,\, \text{and} \quad \delta'^2=-2g.$$

Now consider the $\mathbb Q$-divisor
$$ L_t = tf-\delta'$$
which is ample for $t\gg 0$.  We define
$$ \T = \text{inf}\{ t: L_t \text{ is ample} \}.$$
Clearly $\T\ge \sqrt g$ for if $L_t$ is ample then $0<L_t^2=2t^2-2g$.
In fact conjecturally $\T=\sqrt g$ for ``most'' curves (see Remark
\ref{rmk:generalcurve}). \vspace{.3cm}

Recall that a branched cover from a curve to $\PP^1$ is said to be
simple if it has only ramifications which are locally $z\mapsto z^2$
and no two ramification points map to the same point in $\mathbb P^1$.
Given a simple branched cover $\pi\colon C\to\mathbb P^1$ of degree
$d$, consider the fibre product $C\times_\pi C\subset X$.  This
contains the diagonal $\Delta$ and we let $Z$ be the residual divisor
so $Z+\Delta = C\times_\pi C$.

\begin{lem}\label{lem:classofZ}
  The class of $Z$ in $N^1(X)_\mathbb Q$ is $(d-1)f-\delta'$,  so $Z$
  has self-intersection $ Z^2 = 2(d-1)^2-2g.$
\end{lem}
\begin{proof}
  In the product $\PP^1\times \PP^1$ let $F_i$, $i=1,2$ be the
  numerical class of the coordinate planes and $D$ be the class of the
  diagonal.  Then $D = F_1 + F_2$ so
$$ Z + \delta = C\times_\pi C = (\pi\times\pi)^* (D) = (\pi\times\pi)^*(F_1 + F_2) = df_1+df_2 = df.$$
Hence $Z = df- \delta = (d-1)f-\delta'$ as claimed.
\end{proof}



We will be interested in the case when $2\le d-1<\sqrt g$. Then $Z^2$
is negative, and we will show that $X$ is not slope semistable with
respect to $Z$ for suitable polarisations.

\begin{thm}[Kouvidakis
  \cite{kouvidakis93:_divisors_symmetric}]\label{thm:kouvidakis}
  Suppose $C$ is a smooth curve of genus $g\ge 2$ which admits a
  simple branched cover to $\mathbb P^1$ of degree $d$ with $d-1\le
  \sqrt g $.  Then $\T=\frac{g}{d-1}$.\
\end{thm}
\begin{proof}
See \cite{lazarsfeld04:_positI} Theorem 1.5.8.


\end{proof}


\begin{thm}\label{thm:maintheorem}
  Suppose $C$ is a smooth curve of genus $g$ which admits a
  simple branched cover $\pi\colon C\to \mathbb P^1$ of degree $d$
  with $2\le d-1<\sqrt{g}$.

  Then $X=C\times C$ is not slope semistable with respect $L_t$ for
  $t$ sufficiently close to $\T$.
\end{thm}
\begin{proof}

  Let $t>\T$ so $L_t$ is ample.  The canonical divisor of $X$ is
  $K=(2g-2)f$ so from \eqref{eq:slopeonsurfaces}
  \begin{equation}
    \mu(X,L_t) = -\frac{K.L_t}{L_t^2} = -\frac{t(2g-2)}{t^2-g}.\label{eq:slopeofX}
  \end{equation}

  By \eqref{thm:kouvidakis}, $\T=\frac{g}{d-1}$.  Letting $Z$ be
  the curve from Lemma \ref{lem:classofZ} whose class is
  $(d-1)f-\delta'$ we now bound the Seshadri constant of $Z$.  Since
  $$t-(d-1)>\T-(d-1)= \frac{g}{d-1}-(d-1)>0$$
  we have that $L_t-Z = tf-\delta' - ((d-1)f-\delta') = (t-(d-1))f$ is
  ample.  Thus $\epsilon(Z,L_t)\ge 1$.

  To calculate the slope of $Z$ we need the quantities
\begin{eqnarray} \label{eq:intersectionnumbers}
  L_t.Z &=& (tf-\delta').((d-1)f-\delta') \\
&=& 2t(d-1) -2g, \nonumber \\ 
K.Z&=& (2g-2)f.((d-1)f-\delta')= 2(2g-2)(d-1),  \nonumber \\
Z^2 &=&  ((d-1)f-\delta')^2 =  2(d-1)^2-2g. \nonumber
\end{eqnarray}
Thus from \eqref{eq:slopeonsurfaces},
\begin{eqnarray}
  \mu_1(\O_Z,L_t) &= &\frac{3(2L_t.Z-(K.Z+Z^2))}{2(3L_t.Z-Z^2)} \label{eq:slopeofZ} \\
  &=& \frac{3( 4t(d-1) - 4g - 2(2g-2)(d-1) - 2(d-1)^2+2g)}{2(6t(d-1)-6g-2(d-1)^2+2g)} \nonumber.
\end{eqnarray}
We claim that $\mu_1(\O_Z,L_t)<\mu(X,L_t)$ as $t$ tends to
$\T=\frac{g}{d-1}$ from above.  Since this is an open condition it is
sufficient to show that it holds when $t=\T$.  By (\ref{eq:slopeofX},
\ref{eq:slopeofZ}),
\begin{eqnarray*}
\mu(X,L_{\T}) &=& -\frac{(2g-2)(d-1)}{g-(d-1)^2}, \\
\mu_1(\O_Z,L_{\T})&=&\frac{3(g-(2g-2)(d-1)-(d-1)^2)}{2(g-(d-1)^2)}
\end{eqnarray*}
(notice that our assumption $d-1<\sqrt g$ ensures that both of these
are finite).  Hence as $d-1\ge 2$,
  \begin{align*}
2(g-(d-1)^2)&.(\,\mu_1(\O_Z,L_{\T}) - \mu(X,L_{\T})\,)\\ 
&=3g-3(2g-2)(d-1)-3(d-1)^2 +2(2g-2)(d-1) \\
&=3g-(2g-2)(d-1) -3(d-1)^2\\
&\le 3g-2(2g-2)-12<0.
  \end{align*}
  Thus $\epsilon(Z,L_t)\ge 1$ and
$\mu_1(\O_Z,L_t)<\mu(X,L_t)$ as $t$ tends to $\T$ from above, which
proves that $(X,L_t)$ is not slope semistable.
\end{proof}

\begin{thm}\label{cor:main}
  For $g\ge 5$ there exist smooth curves $C$ of genus $g$ such that
  $X=C\times C$ is not slope semistable with respect to certain
  polarisations.  Thus there are K\"ahler classes on $X$ that do not
  admit K\"ahler metrics of constant scalar curvature.
\end{thm}
\begin{proof}
  By the Riemann existence theorem there exist smooth curves of genus
  $g$ which admit a simple branched covering over $\mathbb P^1$ of
  degree $2\le d-1<\sqrt g$ (in fact one can even take $d=3$).  So by
  \eqref{thm:maintheorem}, $X=C\times C$ is not slope stable
  with respect to certain polarisations and thus not K-semistable by
  \eqref{thm:slopeimpliesKstability}.  The application to cscK
  metrics comes from \eqref{thm:donaldsoncscK}.
\end{proof}

\begin{thm}\label{thm:unstablehilbertchow}
  For $g\ge 5$ there exist smooth curves $C$ of genus $g$ such that
  $X=C\times C$ is not asymptotically Hilbert semistable (resp.\ not
  asymptotically Chow semistable) with respect to suitable
  polarisations.
\end{thm}
\begin{proof}
  This follows from \eqref{cor:main} and
  \eqref{thm:slopeimplieshilbertchowstability}.
\end{proof}

\begin{rmk}
  Let $(X_i,L_i)$ $i=1,2$ be polarised manifolds and $\pi_i\colon
  X_1\times X_2\to X_i$ be the projection maps.  If $L=\pi_1^*L_1+
  \pi_2^*L_2$ then
  $$\mu(X_1\times X_2,L)= \mu(X_1,L_1) + \mu(X_2,L_2).$$
  Moreover if
  $Z$ is subscheme of $X_1$ then one can calculate the slope of
  $Z\times X_2\subset X_1\times X_2$ is $\mu_c(\O_{Z\times X_1}, L) =
  \mu_c(\O_{Z},L_1) + \mu(X_2,L_2)$.  Thus if $(X_1,L_1)$ is slope
  unstable so is the product $(X_1\times X_2,L)$.  
  
  So by taking the product of an slope unstable surface with any
  manifold with negative first Chern class, we get manifolds of any
  dimension $n\ge 2$ with negative first Chern class which have
  K\"ahler classes that do not admit cscK metrics.
\end{rmk}

\section{The J-flow on products of curves}\label{sec:mabuchi}

The Mabuchi functional for a given K\"ahler class $\Omega$ on a
complex manifold $X$ has as its critical points the metrics which are
cscK.  Conjecturally the existence of a cscK metric in $\Omega$ is
equivalent to the properness of the Mabuchi functional.  This is known
to be true when $\Omega$ is proportional to the canonical class
\cite{bandomabuchi87:_uniquen_einst_kaehl,tian97:kaehler_einstein_metrics_positive};
and when $\Omega$ admits a cscK metric the Mabuchi functional is
necessarily bounded from below
\cite{chentian.:_geomet_kaehl,donaldson05:_scalar_ii}

Now suppose that $X$ has negative first Chern class.  Chen \cite{chen00:_mabuchilowerbound} introduces
a flow on K\"ahler manifolds, called the J-flow, and points out that
convergence of this flow implies lower boundedness of the Mabuchi
functional.  In \cite{weinkove:_jflowmabuchi} it is shown that on a
surface with negative first Chern class, the J-flow converges as long
as the class $-2\left(\int_X c_1(X).\Omega\right) \Omega + \left(\int_X \Omega^2\right) c_1(X)$ is
positive.

The following theorem is a particular case of Theorems 1.2 and 1.4 in
\cite{song_winkove:_jflowmabuchi} applied to polarised surfaces.

\begin{thm}\label{thm:jflow}
  Let $(X,L)$ be a polarised surface and define a divisor by
$$\W= 2(K.L) L - (L^2)K.$$
\begin{itemize}
\item If $\W$ is ample then the J-flow converges and the Mabuchi
  functional is proper on the class $c_1(L)$.
\item If $\W$ is not ample then there exist $m$ irreducible curves
  $E_i$ of negative self intersection and positive numbers $a_i$ such
  that
  \begin{equation}
    \W - \sum_{i=1}^m a_iE_i \text{ is ample}. \label{eq:wminusexceptional}
  \end{equation}
\end{itemize}
\end{thm}

In fact if $\W$ is not ample then the J-flow ``blows up'' along the
intersection of all divisors in the linear series $|\Sigma a_iE_i|$
\cite{song_winkove:_jflowmabuchi}.

From this theorem one might expect that $(X,L)$ is stable if $\W$ is
ample, and perhaps that if $\W$ is not ample that the $E_i$ witness
instability of $(X,L)$ (c.f.\ \cite{song_winkove:_jflowmabuchi} Remark
4.7).

When $C$ is a curve, $X=C\times C$ and $L=L_t=tf-\delta'$ with $t>\T$
it is easy to determine when $\W$ is ample.

\begin{lem}\label{lem:w} Let $L=L_t=tf-\delta '$ with $t>\T$.  Then $\W=2(K.L)L-(L^2)K$ is ample if and
  only if $t^2+g>2t\T$ if and only if $t> \T + \sqrt{\T^2-g}$.
\end{lem}
\begin{proof}
  As $L_t^2 = 2t^2-2g$ and $K.L_t=2t(2g-2)$ we have
  \begin{eqnarray*}\label{eq:W}
    \W &=& 4t(2g-2)(tf-\delta') - (2t^2-2g)(2g-2)f \\
&=& 2(2g-2)( (t^2+g)f-2t\delta')
  \end{eqnarray*}
  which is ample if and only if $t^2+g>2t\T$ which occurs if and only
  if $h(t) = t^2-2t\T+g>0$.  But the roots of $h$ are $\T\pm
  \sqrt{\T^2-g}$ and $\T\ge \sqrt g$ so the lemma follows.
\end{proof}

  Thus one can deduce properness of the Mabuchi functional when
  $\T=\sqrt g$.

\begin{cor}
  If $\T=\sqrt{g}$ then the Mabuchi function is proper on any class on
  $X=C\times C$ of the form $c_1(L_t)$ for $t>\T$.
\end{cor}

\begin{rmk}\label{rmk:generalcurve}
  If $g$ is a perfect square and $C$ is very general curve of
  genus $g$ then $\T=\sqrt{g}$ (\cite{lazarsfeld04:_positI} Corollary
  1.5.9).  If the Nagata conjecture holds and $g\ge 10$ then
  the same conclusion holds without the hypothesis that $g$ is a
  perfect square \cite{ciliberto+kouvidakis99_onthesymmetric}.  
\end{rmk}

\begin{rmk}
  Let $C$ be as in \eqref{thm:maintheorem} and $X=C\times C$.
  It would be interesting to know if the Mabuchi functional is bounded
  on the class $c_1(L_t)$ if and only if $t>\T+\sqrt{\T^2-g}$.

  It is not the case that the curve $Z=(d-1)f-\delta'$ from  \eqref{thm:maintheorem} always slope destabilises when
  $t<\T+\sqrt{\T^2-g}$.  For example let $C$ be a curve of genus $g=5$
  admitting a simple branched cover to $\mathbb P^1$ of degree $d=3$.
  Then $\T=5/2$ so $\T+\sqrt{\T^2-g}=\frac{5+\sqrt{5}}{2}>3$.  Put
  $t=3$ so $L_3=3f-\delta'$ is ample.  The slope of $Z$ is
  $$\mu_c(\O_Z,L) = \frac{3(2-15c)}{2c(3+c)}$$
  which is greater than
  $\mu(X,L)= -6$ for all $c>0$.  
  
  By the previous lemma $\alpha$ is not ample.  It is unfortunately not the
  case that any divisor satisfying \eqref{eq:wminusexceptional}
  necessarily destabilise (c.f.\ \cite{song_winkove:_jflowmabuchi}
  Remark 4.7) because in the example above, $\alpha-33Z= 158f-63\delta'$ is
  ample.
\end{rmk}

\bibliography{bibliography}\vspace{5mm}

{\small \noindent {\tt jaross@math.columbia.edu}} \newline
\noindent Department of Mathematics, Columbia University, New York, NY 10027.
USA. \\

\end{document}